\documentclass[english]{article}
\usepackage[T1]{fontenc}
\usepackage[latin9]{inputenc}
\usepackage{mathrsfs}
\usepackage{amsmath}
\usepackage{amssymb}
\usepackage{stackrel}
\makeatletter
\newcommand{\lyxaddress}[1]{
\par {\raggedright #1
\vspace{1.4em}
\noindent\par}
}

\makeatother

\usepackage{babel}
\begin{document}

\title{Soft topology redefined}

\author{Moumita Chiney$^1$, S. K. Samanta$^2$}
\maketitle
\lyxaddress{Department of Mathematics, Visva-Bharati, Santiniketan-731235, West
Bengal, India.
e-mail: $^{1}$ moumi.chiney@gmail.com,$^{2}$ syamal\_123@yahoo.co.in}
\begin{abstract}
In this paper we give a new definition of soft topology using elementary
union and elementary intersection although these operations are not
distributive. Also we have shown that this soft topology is different
from Naz's soft topology and studied some basic properties of this
new type of soft topology. Here we use elementary complement of soft
sets, though law of excluded middle is not valid in general for this
type of complementation. 
\end{abstract}
$\mathbf{Keywords:}$ Soft sets, soft elements, soft topological spaces,
soft base, soft functions, soft continuous functions, soft seperation
axioms.

\section{Introduction}

The concept of soft set was introduced by Molodtsov \cite{key-22}
as a new mathematical tool for dealing with uncertainties. He has
shown several applications of soft set in solving many practical problems
in economics, engineering, social sciences, medical sciences, etc.
Later, Maji et al. defined several operations over soft sets and applied
soft set in decision making problems in \cite{key-18,key-19}. Soft
sets are convenient to be applied in practice and this theory has
potential application in many different fields \cite{key-22} such
as smoothness of functions, game theory, Riemann integration, Perron
integration, probability theory and measure theory. Several authors are trying to develope different mathematical structures in soft set settings (for references please see \cite{key-1,key-2,key-3,key-4,key-5,key-6,key-8,key-9,key-10,key-12,key-13,key-14,key-15,key-16,key-17,key-20,key-23,key-25,key-26,key-27,28,key-29,key-30}). In \cite{key-29}, Shabir and Naz gave a definition of soft topology. Very recently Shi et al. in \cite{key-31}
has commented that soft topology in the sense of Shabir \& Naz \cite{key-29}
can be interpreted as a crisp topology. In this paper, we define soft
topology in different perspective using elementary union and elementary
intersection which is different from that of Naz's topology and study
some basic properties of this new type soft topological space. The
importance of this study lies on the fact that the operations of elementary
union and elementary intersection are not distributive. Also if we
take a soft set and its elementary complement then law of excluded-middle
is not valid in general. This has been mentioned in $Note\:2.24.$\\

\section{Preliminaries}

\textbf{$\mathbf{Definition\,2.1.}$ }\cite{key-22} Let $X$ be a
universal set and $E$ be a set of parameters. Let $P(X)$ denote
the power set of $X$ and $A$ be a subset of $E.$ A pair $(F,A)$
is called a soft set over $X,$ where $F$ is a mapping given by $F:A\rightarrow P(X)$.
In other words, a soft set over $X$ is a parametrized family of subsets
of the universe $X.$ For $\alpha\in A,$ $F(\alpha)$ may be considered
as the set of $\alpha$ - approximate elements of the soft set $(F,A).$$\newline$
In \cite{key-21} the soft sets are redefined as follows:\\
Let $E$ be the set of parameters and $A\subseteq E$. Then for each
soft set $(F,A)$ over $X$ a soft set $(H,E)$ is constructed over
$X$, where $\forall\alpha\in E,$\\ 
$H(\alpha)=\begin{cases}
F(\alpha) & if\:\alpha\in A\\
\phi & if\:\alpha\in E\setminus A.
\end{cases}$\\
Thus the soft sets $(F,A)$ and $(H,E)$ are equivalent to each other
and the usual set operations of the soft sets $(F_{i},A_{i}),i\in\Delta$
is the same as those of the soft sets $(H_{i},E),i\in\Delta$. For
this reason, in this paper, we have considered our soft sets over
the same parameter set $A$. $\newline$\\
Following Molodtsov and Maji et al. \cite{key-18,key-19,key-22} definitions
of soft subset, absolute soft set, null soft set, arbitrary union,
arbitrary intersection of soft sets etc. are presented in \cite{key-26}
considering the same parameter set.$\newline$\\
$\mathbf{Definition\:2.2.}$ \cite{key-26} Let $(F,A)$
and $(G,A)$  be two soft sets over a common universe$X$. \\
$(a)$ $(F,A)$ is said to be a soft subset of $(G,A)$ if $F(\alpha)$
is a subset of $G(\alpha),$ $\forall\alpha\in A.$$\newline$\\
$(b)$  $(F,A)$ and $(G,A)$ are said to be soft equal if $(F,A)$ is a soft subset of $(G,A)$
and $(G,A)$ is a soft subset of $(F,A).$ $\newline$ \\
$(c)$ The complement or relative complement of a soft set $(F,A)$
is denoted by $(F,A)^{C}$ and is defined by $(F,A)^{C}=(F^{C},A),$ where $F^{C}(\alpha)=X\setminus F(\alpha),\forall\alpha\in A.$$\newline$\\
$(d)$ (Null soft set)  $(F,A)$ over $X$ is said to be
null soft set denoted by $(\tilde{\Phi},A)$, where $F(\alpha)=\phi,\forall\alpha\in A.$$\newline$\\
$(e)$ (Absolute soft set) $(F,A)$ over $X$ is said to
be absolute soft set denoted by $(\tilde{X},A),$ if $F(\alpha)=X,,\forall\alpha\in A.$
$\newline$\\
 $(f)$ Union of  $(F,A)$ and $(G,A)$  is denoted by $(F,A)\tilde{\cup}(G,A)$ and defined by
$\left[(F,A)\tilde{\cup}(G,A)\right](\alpha)=F(\alpha)\cup G(\alpha),\forall\alpha\in A.$$\newline$\\
 $(g)$ Intersection of $(F,A)$ and $(G,A)$  is denoted by $(F,A)\tilde{\cap}(G,A)$ and defined
by $\left[(F,A)\tilde{\cap}(G,A)\right](\alpha)=F(\alpha)\cap G(\alpha),\forall\alpha\in A.$
$\newline$\\
$\mathbf{Definition\:2.3.}$ \cite{key-29} Let $\tau$ be a collection
of soft sets over $X$. Then $\tau$ is said to be a soft topology
on $X$ if\\
$(i)$ $(\tilde{\Phi},A),$ $(\tilde{X},A)\in\tau$. \\
$(ii)$ the intersection of any two soft sets in $\tau$ belongs to
$\tau$. \\
$(iii)$ the union of any number of soft sets in $\tau$ belongs to
$\tau$. 
\noindent The triplet $(X,A,\tau)$ is called a soft topological space
over $X$. $\newline$\\
$\mathbf{Definition\:2.4.}$ \cite{key-13} Let $\tau$ be a collection
of soft sets over $X$. Define $\tau(\alpha)=\{F(\alpha):(F,A)\in\tau\}$
for $\alpha\in A.$ Then $\tau$ is said to be a topology of soft
subsets over $(X,A)$ if $\tau(\alpha)$ is a crisp topology on $X,$
$\forall\alpha\in A.$ In this case $((X,A),\tau)$ is said to be
a topological space of soft subsets. If $\tau$ is a topology of soft
susbsets over $(X,A)$, then the members of $\tau$ are called open
soft sets and a soft set $(F,A)$ over $X$ is said to be soft closed
if $(F,A)^{C}\in\tau.$$\newline$\\
\textbf{$\mathbf{Definition\:2.5.}$ }\cite{key-7} Let $X$ be a
non-empty set and $A$ be a non-empty parameter set. Then a function
$\tilde{x}:A\rightarrow X$ is said to be a soft element of $X.$
A soft element $\tilde{x}$ of $X$ is said to belong to a soft set
$(F,A)$ of $X$, which is denoted by $\tilde{x}\tilde{\in}(F,A),$
if $\tilde{x}(\lambda)\in F(\lambda),$ for all $\lambda\in A.$ Thus
for a soft set $(F,A)$ over $X$ with respect to the index set $A$
with $F(\lambda)\neq\phi,$ for all $\lambda\in A$, we have $F(\lambda)=\{\tilde{x}(\lambda):\tilde{x}\tilde{\in}(F,A)\},$
for all $\lambda\in A$. \\
By the notation $\bar{x}$ we will denote a particular type of soft
element such that $\bar{x}(\lambda)=x,$ for all $\lambda\in A.$
$\newline$\\
Let $X$ be an initial universal set and $A$ be a non-empty parameter
set. Throughout the paper we consider the null soft set $(\tilde{\Phi},A)$
and those soft sets $(F,A)$ over $X$ for which $F(\alpha)\neq\phi,\forall\alpha\in A.$
We denote this collection by $S(\tilde{X}).$ Thus for $(F,A)[\neq(\tilde{\Phi},A)]\in S(\tilde{X}),F(\alpha)\neq\phi,\forall\alpha\in A.$$\newline$\\
$\mathbf{Proposition\:2.6.}$ \cite{key-9} Any collection of soft
elements of a soft set can generate a soft subset of that soft set.\\
The soft set constructed from a collection of soft elements \ensuremath{\mathscr{B}}
will be denoted by $SS(\text{\ensuremath{\mathscr{B}}}).$ For any
soft set $(F,A)\in S(\tilde{X}),$ the collection of all soft elements
of $(F,A)$ is denoted by $SE(F,A).$$\newline$\\
$\mathbf{Proposition\:2.7.}$ \cite{key-9} For any soft sets $(F,A),(G,A)\in S(\tilde{X}),$
$(F,A)\tilde{\subseteq}(G,A)$ iff every soft element of $(F,A)$
is also a soft element of $(G,A).$$\newline$\\
$\mathbf{Definition\:2.8.}$ \cite{key-9} For any two soft sets
$(F,A),(G,A)\in S(\tilde{X}),$ \\
$(i)$ elementary union of $(F,A)$ and $(G,A)$ is denoted by $(F,A)\Cup(G,A)$
and is defined by $(F,A)\Cup(G,A)=SS(\text{\ensuremath{\mathscr{B}}}),$
where $\text{\ensuremath{\mathscr{B}}}=\{\tilde{x}\tilde{\in}(\tilde{X},A):\tilde{x}\tilde{\in}(F,A)\; or\;\tilde{x}\tilde{\in}(G,A)\};$
i.e. $(F,A)\Cup(G,A)=SS(SE(F,A)\cup SE(G,A)).$\\
$(ii)$ elementary intersection of $(F,A)$ and $(G,A)$ is denoted
by $(F,A)\Cap(G,A)$ and is defined by $(F,A)\Cap(G,A)=SS(\text{\ensuremath{\mathscr{B}}}),$
where \\
$\text{\ensuremath{\mathscr{B}}}=\{\tilde{x}\tilde{\in}(\tilde{X},A):\tilde{x}\tilde{\in}(F,A)\; and\;\tilde{x}\tilde{\in}(G,A)\};$
i.e. $(F,A)\Cap(G,A)=SS(SE(F,A)\cap SE(G,A)).$$\newline$\\
$\mathbf{Definition\:2.9.}$ \cite{key-9} For any soft set $(F,A)\in S(\tilde{X}),$
the elementary complement of $(F,A)$ is denoted by $(F,A)^{\mathbb{C}}$
and is defined by $(F,A)^{\mathbb{C}}=SS(\text{\ensuremath{\mathscr{B}}}),$
where $\text{\ensuremath{\mathscr{B}}}=\{\tilde{x}\tilde{\in}(\tilde{X},A):\tilde{x}\tilde{\in}(F,A)^{C}\}$
and $(F,A)^{C}$ is the complement of $(F,A).$$\newline$\\
$\mathbf{Remark\:2.10.}$ \cite{key-9} It can be easily verified
that if $(F,A),(G,A)\in S(\tilde{X}),$ then $(F,A)\Cup(G,A)$, $(F,A)\Cap(G,A)$
and $(F,A)^{\mathbb{C}}$ are members of $S(\tilde{X}).$$\newline$\\
$\mathbf{Proposition\:2.11.}$ \cite{key-9} For any two soft sets
$(F,A),(G,A)\in S(\tilde{X}),$ \\
$(i)$ $(F,A)\Cup(G,A)=(F,A)\tilde{\cup}(G,A).$\\
$(ii)$ $(F,A)\Cap(G,A)=(F,A)\tilde{\cap}(G,A),$ if $(F,A)\Cap(G,A)\neq(\tilde{\Phi},A).$$\newline$\\
$\mathbf{Remark\:2.12.}$ The above results can be extended easily
to arbitrary union and arbitrary intersection. $\newline$ \\
$\mathbf{Proposition\:2.13.}$ \cite{key-9} For any soft set $(F,A)\in S(\tilde{X}),$
in general, $(F,A)^{\mathbb{C}}\tilde{\subseteq}(F,A)^{C}$ and $(F,A)^{\mathbb{C}}=(F,A)^{C},$
if $(F,A)^{\mathbb{C}}\neq(\tilde{\Phi},A)$ i.e. if $(F,A)^{C}\in S(\tilde{X}).$
$\newline$\\
$\mathbf{Proposition\:2.14.}$ \cite{key-9} For any soft set $(F,A)\in S(\tilde{X}),$
\\
$(i)$ $(F,A)\Cap(F,A)^{\mathbb{C}}=(\tilde{\Phi},A)$,\\
$(ii)$ In general,$(F,A)\Cup(F,A)^{\mathbb{C}}\tilde{\subseteq}(\tilde{X},A)$
but if $(F,A)^{\mathbb{C}}\neq(\tilde{\Phi},A)$, then $(F,A)\Cup(F,A)^{\mathbb{C}}=(\tilde{X},A).$
$\newline$\\
$\mathbf{Proposition\:2.15.}$ \cite{key-9} Let $\{(Y_{i},A),i\in\varLambda\}$
be a family of soft sets and $\{\text{\ensuremath{\mathscr{B}}}_{i},i\in\varLambda\}$
be a family of collection of soft elements such that $(Y_{i},A)=SS(\text{\ensuremath{\mathscr{B}}}_{i}),\forall i\in\varLambda.$
Then $\underset{i\in\varLambda}{\Cup}(Y_{i},A)=SS(\underset{i\in\varLambda}{\cup}\text{\ensuremath{\mathscr{B}}}_{i})$.
$\newline$\\
$\mathbf{Proposition\:2.16.}$ Let $\{(Y_{i},A),i\in\varLambda\}$
be a family of soft sets. Then $\underset{i\in\varLambda}{\Cap}(Y_{i},A)=\underset{i\in\varLambda}{\tilde{\cap}}(Y_{i},A)$
iff $\underset{i\in\varLambda}{\tilde{\cap}}(Y_{i},A)\in S(\tilde{X}).$
$\newline$\\
$\mathbf{Proposition\:2.17.}$ Let $\{(Y_{i},A),i\in\varLambda\}$
be a family of soft sets and $\{\text{\ensuremath{\mathscr{B}}}_{i},i\in\varLambda\}$
be a family of collection of soft elements such that $(Y_{i},A)=SS(\text{\ensuremath{\mathscr{B}}}_{i}),\forall i\in\varLambda.$
Then $\underset{i\in\varLambda}{\Cap}(Y_{i},A)\tilde{\supseteq}SS(\underset{i\in\varLambda}{\cap}\text{\ensuremath{\mathscr{B}}}_{i})$.
$\newline$\\
$Proof.$ We have $\text{\ensuremath{\mathscr{B}}}_{i}\subset SE(Y_{i},A),\forall i\in\varLambda$\\
$\Rightarrow\underset{i\in\varLambda}{\cap}\text{\ensuremath{\mathscr{B}}}_{i}\subset\underset{i\in\varLambda}{\cap}SE(Y_{i},A)$\\
$\Rightarrow SS(\underset{i\in\varLambda}{\cap}\text{\ensuremath{\mathscr{B}}}_{i})\tilde{\subset}SS[\underset{i\in\varLambda}{\cap}SE(Y_{i},A)]$\\
$\Rightarrow SS(\underset{i\in\varLambda}{\cap}\text{\ensuremath{\mathscr{B}}}_{i})\tilde{\subset}\underset{i\in\varLambda}{\Cap}(Y_{i},A)$.
$\newline$\\
$\mathbf{Example\:2.18.}$ In general, $\underset{i\in\varLambda}{\Cap}(Y_{i},A)\neq SS(\underset{i\in\varLambda}{\cap}\text{\ensuremath{\mathscr{B}}}_{i})$.
$\newline$\\
For example, let $X=\{x,y\},A=\{\lambda,\mu\},\text{\ensuremath{\mathscr{B}}}_{1}=\{\tilde{x}_{1},\tilde{x}_{2},\tilde{x}_{3}\},\text{\ensuremath{\mathscr{B}}}_{2}=\{\tilde{x}_{4}\};$
where $\tilde{x}_{1}(\lambda)=\{x\},$ $\tilde{x}_{1}(\mu)=\{y\};$
$\tilde{x}_{2}(\lambda)=\{x\},$ $\tilde{x}_{2}(\mu)=\{x\};$ $\tilde{x}_{3}(\lambda)=\{y\},$
$\tilde{x}_{3}(\mu)=\{y\};$$\tilde{x}_{4}(\lambda)=\{y\},$ $\tilde{x}_{4}(\mu)=\{x\}.$
Also, let $(Y_{1},A)=SS\{\text{\ensuremath{\mathscr{B}}}_{1}\},(Y_{2},A)=SS\{\text{\ensuremath{\mathscr{B}}}_{2}\}$.
Then $\text{\ensuremath{\mathscr{B}}}_{1}\cap\text{\ensuremath{\mathscr{B}}}_{2}=\phi\: i.e.\: SS(\text{\ensuremath{\mathscr{B}}}_{1}\cap\text{\ensuremath{\mathscr{B}}}_{2})=(\tilde{\varPhi},A)$.
But $(Y_{1},A)\Cap(Y_{2},A)=(Y_{2},A).$$\newline$\\
$\mathbf{Proposition\:2.19.}$ Let $(Y,A),(Z,A)\in S(\tilde{X}).$
Then \\
$(i)$ $SE[(Y,A)\Cap(Z,A)]=SE(Y,A)\cap SE(Z,A)$.\\
$(ii)$ $SE[(Y,A)\Cup(Z,A)]\supset SE(Y,A)\cup SE(Z,A)$. $\newline$\\
$\mathbf{Example\:2.20.}$ In general, $SE[(Y,A)\Cup(Z,A)]\neq SE(Y,A)\cup SE(Z,A)$.
\\
 Consider $X=\{x,y,z\},$ $A=\{\lambda,\mu\}$, $(Y,A),(Z,A)\in S(\tilde{X}),$
where $Y(\lambda)=\{y\},$ $Y(\mu)=\{z\};$ $Z(\lambda)=\{z\},$ $Z(\mu)=\{x\}.$
Also consider the soft elements $\tilde{x}_{1},$ $\tilde{x}_{2},$
$\tilde{x}_{3},$ $\tilde{x}_{4},$ where $\tilde{x}_{1}(\lambda)=\{y\},$
$\tilde{x}_{1}(\mu)=\{z\};$ $\tilde{x}_{2}(\lambda)=\{z\},$ $\tilde{x}_{2}(\mu)=\{x\};$
$\tilde{x}_{3}(\lambda)=\{y\},$ $\tilde{x}_{3}(\mu)=\{x\};$$\tilde{x}_{4}(\lambda)=\{z\},$
$\tilde{x}_{4}(\mu)=\{z\}.$ Then $SE(Y,A)\cup SE(Z,A)=\{\tilde{x}_{1},\tilde{x}_{2}\}$,
but $SE[(Y,A)\Cup(Z,A)]=\{\tilde{x}_{1},\tilde{x}_{2},\tilde{x}_{3},\tilde{x}_{4}\}$.$\newline$\\
$\mathbf{Remark\:2.21.}$ The operations of elementary union and elementary
intersection are not distributive over $S(\tilde{X}).$ This is shown
in the following example.$\newline$\\
$\mathbf{Example\:2.22.}$ Let $X=\{x,y,z\}$ and $A=\{\alpha,\beta\}$.
Now consider three soft sets $(F,A),$ $(G,A)$ and $(H,A)\in S(\tilde{X})$,
where $F(\alpha)=\{x\},$ $F(\beta)=\{y\};$ $G(\alpha)=\{y\},$ $G(\beta)=\{z\}$
and $H(\alpha)=\{y\},$ $H(\beta)=\{y\}.$ \\
Then $[(F,A)\Cup(G,A)]\Cap(H,A)\neq[(F,A)\Cap(H,A)]\Cup[(G,A)\Cap(H,A)]$
as\\
$(F,A)\Cap(H,A)=(\tilde{\Phi},A)$, $(G,A)\Cap(H,A)=(\tilde{\Phi},A)$
i.e. $[(F,A)\Cap(H,A)]\Cup[(G,A)\Cap(H,A)]=(\tilde{\Phi},A)$ but
$[(F,A)\Cup(G,A)]\Cap(H,A)=(H,A)$.$\newline$\\
Also $[(F,A)\Cap(H,A)]\Cup(G,A)\neq[(F,A)\Cup(G,A)]\Cap[(H,A)\Cup(G,A)]$
as $[(F,A)\Cup(G,A)](\alpha)=\{x,y\},$ $[(F,A)\Cup(G,A)](\beta)=\{y,z\}$
and $[(H,A)\Cup(G,A)](\alpha)=\{y\},$ $[(H,A)\Cup(G,A)](\beta)=\{y,z\}$
i.e. $[(F,A)\Cup(G,A)]\Cap[(H,A)\Cup(G,A)](\alpha)=\{y\},$ $[(F,A)\Cup(G,A)]\Cap[(H,A)\Cup(G,A)](\beta)=\{y,z\}$
but $[(F,A)\Cap(H,A)]\Cup(G,A)=(G,A).$$\newline$\\
$\mathbf{Proposition\:2.23.}$ $(i)$ If $(F,A),(F,A)^{C}\in S(\tilde{X}),$
then $(F,A)^{C}=(F,A)^{\mathbb{C}}$ i.e. then $[(F,A)^{\mathbb{C}}]^{\mathbb{C}}=(F,A)$.\\
$(ii)$ Let $\{(F_{i},A):i\in\triangle\}$ be any collection of soft
sets in $S(\tilde{X}),$ then 
\begin{itemize}
\item[(a)] $\left[\Cup\{(F_{i},A):i\in\triangle\}\right]^{\mathbb{C}}=\Cap\left\{ (F_{i},A)^{\mathbb{C}}:i\in\triangle\right\}.$
\item[(b)] $\left[\Cap\{(F_{i},A):i\in\triangle\}\right]^{\mathbb{C}}=\Cup\left\{ (F_{i},A)^{\mathbb{C}}:i\in\triangle\right\}.$
\end{itemize}
$\mathbf{Note\:2.24.}$ Following observations are to be noted:
\begin{itemize}
\item[(i)] The operations of elementary union and
elementary intersection are not distributive over $S(\tilde{X}).$
\item[(ii)] In general,$(F,A)\Cup(F,A)^{\mathbb{C}}\tilde{\subseteq}(\tilde{X},A)$,
$(F,A)\Cup(F,A)^{\mathbb{C}}\neq(\tilde{X},A).$ However if $(F,A)^{\mathbb{C}}\neq(\tilde{\Phi},A)$,
then $(F,A)\Cup(F,A)^{\mathbb{C}}=(\tilde{X},A).$
\item[(iii)] In general, $\left((F,A)^{\mathbb{C}}\right)^{\mathbb{C}}\neq(F,A).$
However if $(F,A)^{\mathbb{C}}\neq(\tilde{\Phi},A)$, then $\left((F,A)^{\mathbb{C}}\right)^{\mathbb{C}}$
$=(F,A).$
\end{itemize}
\section{Soft topology}
\noindent $\mathbf{Definition\:3.1.}$ Let $\tau$ be the collection
of soft sets of $S(\tilde{X}).$ Then $\tau$ is said to be a soft
topology on $(\tilde{X},A)$ if \\
$(i)$ $(\tilde{\Phi},A)$ and $(\tilde{X},A)$ belong to $\tau$.\\
$(ii)$ the elementary union of any number of soft sets in $\tau$
belongs to $\tau.$\\
$(iii)$ the elementary intersection of two soft sets in $\tau$ belongs
to $\tau.$ 
\noindent The triplet $(\tilde{X},\tau,A)$ is called a soft topological
space.$\newline$\\
$\mathbf{Definition\:3.2.}$ In a soft topological space$(\tilde{X},\tau,A),$
the members of $\tau$ are called soft open sets in $(\tilde{X},\tau,A).$$\newline$\\
$\mathbf{Definition\:3.3.}$ Let $(\tilde{X},\tau,A)$ be a soft topological
space. A soft set $(F,A)\in S(\tilde{X})$ is said to be a soft closed
set in $(\tilde{X},\tau,A)$ if its relative complement $(F,A)^{C}\in S(\tilde{X})$
and $(F,A)^{\mathbb{C}}\in\tau.$$\newline$\\
$\mathbf{Proposition\:3.4.}$ Let $(\tilde{X},\tau,A)$ be a soft
topological space. Then \\
$(i)$ $(\tilde{\Phi},A)$ and $(\tilde{X},A)$ are soft closed soft
sets in $(\tilde{X},\tau,A)$.\\
$(ii)$ arbitrary elementary intersection of soft sets is soft closed.$\newline$\\
$Proof.$ $(i)$ Since $(\tilde{\Phi},A)^{C}=(\tilde{X},A)\in S(\tilde{X})$
and $(\tilde{\Phi},A)^{\mathbb{C}}=(\tilde{X},A)\in\tau$ it follows
that $(\tilde{\Phi},A)$ is soft closed.\\
Similarly, $(\tilde{X},A)$ is soft closed. $\newline$ \\
$(ii)$ Let $\{(Y_{i},A):i\in\triangle\}$ be a family of soft closed
sets in $(\tilde{X},\tau,A)$. We have to show that $\Cap\{(Y_{i},A):i\in\triangle\}=(Y,A)$
(say) is soft closed in $(\tilde{X},\tau,A)$. From $Remark\:2.10,$
$(Y,A)\in S(\tilde{X})$\\
If $(Y,A)=(\tilde{\Phi},A),$ then $(Y,A)$ is soft closed in $(\tilde{X},\tau,A)$.
\\
If $(Y,A)\neq(\tilde{\Phi},A),$ then from $Remark\:2.12,$ $\Cap\{(Y_{i},A):i\in\triangle\}=\tilde{\cap}\{(Y_{i},A):i\in\triangle\}$
and hence $(Y,A)^{C}=[\tilde{\cap}\{(Y_{i},A):i\in\triangle\}]^{C}=\tilde{\cup}\{(Y_{i},A)^{C}:i\in\triangle\}$
$=\Cup\{(Y_{i},A)^{C}:i\in\triangle\}\in S(\tilde{X}),$ since $(Y_{i},A)$
is soft closed $(Y_{i},A)^{C}\in S(\tilde{X}),\forall i\in\triangle$.\\
 Now $(Y,A)^{\mathbb{C}}=\left[\Cap\{(Y_{i},A):i\in\triangle\}\right]^{\mathbb{C}}=\Cup\{(Y_{i},A)^{\mathbb{C}}:i\in\triangle\}\in\tau,$
as $(Y_{i},A)^{\mathbb{C}}\in\tau,\forall i\in\triangle$ and $\tau$
is closed under arbitrary elementary union. \\
Therefore $(Y,A)$ is soft closed in $(\tilde{X},\tau,A)$. $\newline$
\\
$\mathbf{Remark\:3.5.}$ Elementary union of two soft closed sets
is not soft closed in general. This is shown in the following example.
$\newline$ \\
$\mathbf{Example\:3.6.}$ Let $X=\{x,y,z\}$ and $A=\{\alpha,\beta\}.$
Then $\tau=\{(\tilde{\Phi},A),(\tilde{X},A),$ $(F,A),(G,A)\}$ where
$F(\alpha)=\{y,z\},F(\beta)=\{x,z\}$ and $G(\alpha)=\{x\},G(\beta)=\{y,z\}.$
Then $\tau$ is a soft topology as per $Definition\;3.1$. The soft
closed sets are $(\tilde{\Phi},A),(\tilde{X},A),$ $(F,A)^{\mathbb{C}},(G,A)^{\mathbb{C}}$,
where $F^{\mathbb{C}}(\alpha)=\{x\},F^{\mathbb{C}}(\beta)=\{y\}$
and $G^{\mathbb{C}}(\alpha)=\{y,z\},G^{\mathbb{C}}(\beta)=\{x\}.$
But $(F,A)^{\mathbb{C}}\Cup(G,A)^{\mathbb{C}}$ is not soft closed.
$\newline$ \\
$\mathbf{Proposition\:3.7.}$ Let $\{(Y_{i},A):i=1,2,....,n\}$ be
a finite family of soft closed sets in $(\tilde{X},\tau,A)$. Then
$\stackrel[i=1]{n}{\Cup}(Y_{i},A)$ is soft closed in $(\tilde{X},\tau,A)$
if $\stackrel[i=1]{n}{\Cap}(Y_{i},A)^{\mathbb{C}}(\neq(\tilde{\Phi},A))\in S(\tilde{X}).$$\newline$\\
$Proof.$ Let $\{(Y_{i},A):i=1,2,....,n\}$ be a finite family of
soft closed sets in $(\tilde{X},\tau,A)$. We have to show that $\stackrel[i=1]{n}{\Cup}(Y_{i},A)=(Y,A)$
is soft closed in $(\tilde{X},\tau,A)$. Since for each $i=1,2,....,n,$
$(Y_{i},A)$ is closed in $(\tilde{X},\tau,A)$, $(Y_{i},A),(Y_{i},A)^{C}\in S(\tilde{X})$
and $(Y_{i},A)^{\mathbb{C}}\in\tau$. \\
Now $(Y,A)^{C}=[\Cup\{(Y_{i},A):i=1,2,..,n\}]^{C}=[\tilde{\cup}\{(Y_{i},A):i=1,2,..,n\}]^{C}=\tilde{\cap}\{(Y_{i},A)^{C}:i=1,2,..,n\}=\Cap\{(Y_{i},A)^{C}:i=1,2,..,n\}\in S(\tilde{X})$,
as $\stackrel[i=1]{n}{\Cap}(Y_{i},A)^{\mathbb{C}}\neq(\tilde{\Phi},A).$\\
Also $(Y,A)^{\mathbb{C}}=[\Cup\{(Y_{i},A):i=1,2,..,n\}]^{\mathbb{C}}=\Cap\{(Y_{i},A)^{\mathbb{C}}:i=1,2,..,n\}$.
Then $(Y,A)^{\mathbb{C}}\in\tau,$ as $(Y_{i},A)^{\mathbb{C}}\in\tau,\forall i=1,2,..,n$
and $\tau$ is closed under finite elementary intersection. Hence
$(Y,A)$ is soft closed in $(\tilde{X},\tau,A)$. $\newline$\\
$\mathbf{Remark\:3.8.}$ $(a)$ Soft topology due to Shabir and Naz
and the soft topology as given by us in $Definition\:3.1$ are different.\\
Let $X=\{x,y,z\}$ and $A=\{\alpha,\beta\}.$ Then $\tau=\{(\tilde{\Phi},A),(\tilde{X},A),(F,A),(G,A)\}$
where $F(\alpha)=\{x,y\},F(\beta)=\{x,z\}$ and $G(\alpha)=\{z\},G(\beta)=\{y,z\}$.
Then $(F,A)\Cup(G,A)=(\tilde{X},A)$ and $(F,A)\Cap(G,A)=(\tilde{\Phi},A).$
Thus $\tau$ is a soft topology as per $Definition\:3.1$. Now as
$(F,A)\tilde{\cap}(G,A)=\{\phi,\{z\}\},$ and $\{\phi,\{z\}\}$ does
not belongs to $\tau$, it is not a soft topology as per Shabir and
Naz. Now if $\tau^{\prime}=\{(\tilde{\Phi},A),(\tilde{X},A),(F,A),(G,A),(H,A)\}$
where $H(\alpha)=\phi,H(\beta)=\{z\}$. Then $\tau^{\prime}$ is a
soft topology as per Shabir and Naz but not a soft topology as per
$Definition\:3.1$, since $(H,A)\notin S(\tilde{X}).$ $\newline$\\
$(b)$ Soft topology due to Hazra et al. and the soft topology as
per $Definition\;3.1$ are different.\\
Consider the universal set $X,$ the parameter set $A$ and $\tau,\tau^{\prime}$
are similar as above. Then $\tau$ is a soft topology as per $Definition\:3.1$.
Here $\tau(\alpha)=\{\phi,X,\{x,y\},\{z\}\}$ and $\tau(\beta)=\{\phi,X,\{x,z\},\{y,z\}\}$.
So, $\tau(\beta)$ is not a crisp topology on $X.$ Thus $\tau$ is
not a topology of soft subsets over $(X,A)$ as per Hazra et al. But
$\tau^{\prime}(\alpha)=\{\phi,X,\{x,y\},\{z\}\}$ and $\tau^{\prime}(\beta)=\{\phi,X,\{x,z\},\{y,z\},\{z\}\}$
are crisp topologies on $X.$ Thus $\tau^{\prime}$ is a topology
of soft subsets over $(X,A)$ as per Hazra et al. but not a soft topology
as per $Definition\:3.1$, as we have seen before.$\newline$\\
However we have the following result:\\
$\mathbf{Proposition\:3.9.}$ Let $(\tilde{X},\tau,A)$ be a soft
topological space as per $Definition\:3.1$ and if $(F,A),(G,A)\in\tau$
$\Rightarrow(F,A)\tilde{\cap}(G,A)\in S(\tilde{X}),$ then $\tau$
defines a soft topology as per Hazra et al. Further if $\tau$ is
a soft topology as per Hazra et al., then $\tau^{\prime}=\{(F,A)\in S(\tilde{X}):F(\alpha)\in\tau(\alpha),\forall\alpha\in A\}$
is a soft topology on $(\tilde{X},A)$ as per $Definition\:3.1$.
$\newline$ \\
$Proof.$ Let $(F,A),(G,A)\in\tau$ and $(F,A)\tilde{\cap}(G,A)\in S(\tilde{X}).$
Then we have to show that $\tau(\alpha)=\{F(\alpha):(F,A)\in S(\tilde{X})\}$
is a crisp topology on $X$. We only show the intersection property
as rest are obvious. Let $F(\alpha),G(\alpha)\in\tau(\alpha)$ for
some $(F,A),(G,A)\in\tau$. \\
Now $F(\alpha)\cap G(\alpha)=\left[(F,A)\tilde{\cap}(G,A)\right](\alpha)$\\
$=\left[(F,A)\Cap(G,A)\right](\alpha)$ {[}$(F,A)\tilde{\cap}(G,A)\in S(\tilde{X})\Rightarrow(F,A)\tilde{\cap}(G,A)=(F,A)\Cap(G,A)${]}\\
$\in\tau(\alpha)$ {[}as $(F,A)\Cap(G,A)\in\tau${]} \\
Thus $\tau$ defines a soft topology as per Hazra et al.\\
Next, Obviously $(\tilde{\Phi},A),(\tilde{X},A)\in\tau^{\prime}.$\\
Now let $(F_{1},A),(F_{2},A)\in\tau^{\prime}.$ If $(F_{1},A)\Cap(F_{2},A)=(\tilde{\Phi},A).$
Then obviously $(F_{1},A)\Cap(F_{2},A)\in\tau^{\prime}.$ \\
Again if $(F_{1},A)\Cap(F_{2},A)\neq(\tilde{\Phi},A).$ Then $(F_{1},A)\Cap(F_{2},A)=(F_{1},A)\tilde{\cap}(F_{2},A)$.
\\
Also $F_{1}(\alpha),F_{2}(\alpha)\in\tau(\alpha),\forall\alpha\in A.$
\\
Thus $\left[(F_{1},A)\Cap(F_{2},A)\right](\alpha)$ $=\left[(F_{1},A)\tilde{\cap}(F_{2},A)\right](\alpha)$
$=F_{1}(\alpha)\cap F_{2}(\alpha)$ $\in\tau(\alpha),\forall\alpha\in A.$\\
Hence, $(F_{1},A)\Cap(F_{2},A)\in\tau^{\prime}.$\\
Similarly if we take $(F_{i},A)\in\tau^{\prime},\forall i\in\triangle,$
then $\Cup\{(F_{i},A),i\in\triangle\}\in\tau^{\prime}.$\\
Therefore, $\tau$ is a soft topology on $(\tilde{X},A)$.$\newline$\\
$\mathbf{Proposition\:3.10.}$ Let $(\tilde{X},\tau_{1},A)$ and $(\tilde{X},\tau_{2},A)$
be two soft topological spaces. Then $(\tilde{X},\tau_{1}\cap\tau_{2},A)$
is a soft topological space.$\newline$\\
\noindent $\mathbf{Remark\:3.11:}$ The union of two soft topologies
on $(\tilde{X},A)$ may not be a soft topology on $(\tilde{X},A)$.
$\newline$\\
$\mathbf{Example\:3.12:}$ Let $X=\{x,y,z\}$ and $A=\{\alpha,\beta\}.$\\
 Let $\tau_{1}=\{(\tilde{\Phi},A),(\tilde{X},A),$ $(F_{1},A),(G_{1},A)\}$
where $F_{1}(\alpha)=\{x,y\},F_{1}(\beta)=\{x,z\}$ and $G_{1}(\alpha)=\{z\},G_{1}(\beta)=\{y,z\}$
and\\
 $\tau_{2}=\{(\tilde{\Phi},A),(\tilde{X},A),(F_{2},A),(G_{2},A)\}$
where $F_{2}(\alpha)=\{y\},F_{2}(\beta)=\{y\}$ and $G_{2}(\alpha)=\{x,z\},G_{2}(\beta)=\{x,z\}.$
\\
Now, we define $\tau=\tau_{1}\cup\tau_{2}$ $=\{(\tilde{\Phi},A),(\tilde{X},A),(F_{1},A),(G_{1},A),(F_{2},A),(G_{2},A)\}.$
\\
If we take $(F_{1},A)\Cup(F_{2},A)=(H,A)$ then $H(\alpha)=\{x,y\},H(\beta)=\{x,y,z\}$.
But $(H,A)\notin\tau.$ \\
Thus $\tau$ is not a soft topology as in $Definition\:3.1.$ $\newline$\\
$\mathbf{Definition\:3.13.}$ Let $(\tilde{X},\tau,A)$ be a soft
topological space and $(F,A)\in S(\tilde{X}).$ Then the soft closure
of $(F,A),$ denoted by $\overline{(F,A)}$ is defined as the elementary
intersection of all soft closed super sets of $(F,A)$. Clearly $\overline{(F,A)}$
is the smallest soft closed set in $(\tilde{X},\tau,A)$ which contains
$(F,A)$ by $Proposition\:3.4.$ $\newline$\\
$\mathbf{Proposition\:3.14.}$ Let $(\tilde{X},\tau,A)$ be a soft
topological space and $(F,A),(G,A)\in S(\tilde{X}).$ Then\\
$(i)$ $\overline{(\tilde{\Phi},A)}=(\tilde{\Phi},A),\overline{(\tilde{X},A)}=(\tilde{X},A).$\\
$(ii)$ $(F,A)\tilde{\subseteq}\overline{(F,A)}$.\\
$(iii)$ $(F,A)$ is soft closed if and only if $(F,A)=\overline{(F,A)}.$\\
$(iv)$ $\overline{\overline{(F,A)}}=\overline{(F,A)}$\\
$(v)$ $(F,A)\tilde{\subseteq}(G,A)\Rightarrow\overline{(F,A)}\tilde{\subseteq}\overline{(G,A)}.$\\
$(vi)$ $\overline{(F,A)}\Cup\overline{(G,A)}\tilde{\subseteq}\overline{(F,A)\Cup(G,A)}.$
Equality holds if, $\left(\overline{(F,A)}\Cup\overline{(G,A)}\right)^{C}\in S(\tilde{X}).$\\
$(vii)$ $\overline{(F,A)\Cap(G,A)}\tilde{\subseteq}\overline{(F,A)}\Cap\overline{(G,A)}$.$\newline$\\
$Proof.$ Proofs of $(i)$-$(v)$ are obvious.$\newline$\\
$(vi)$ Since $(F,A)\tilde{\subseteq}(F,A)\Cup(G,A)$ and $(G,A)\tilde{\subseteq}(F,A)\Cup(G,A),$
so by $(v),$\\
$\overline{(F,A)}\tilde{\subseteq}\overline{(F,A)\Cup(G,A)}$ and
$\overline{(G,A)}\tilde{\subseteq}\overline{(F,A)\Cup(G,A)}$. \\
Thus $\overline{(F,A)}\Cup\overline{(G,A)}\tilde{\subseteq}\overline{(F,A)\Cup(G,A)}.$
\\
Let, $\left(\overline{(F,A)}\Cup\overline{(G,A)}\right)^{C}\in S(\tilde{X}).$
Now,$(F,A)\tilde{\subseteq}\overline{(F,A)}$ and $(G,A)\tilde{\subseteq}\overline{(G,A)}$.
Then $(F,A)\Cup(G,A)\tilde{\subseteq}\overline{(F,A)}\Cup\overline{(G,A)}.$
Also, $\left(\overline{(F,A)}\Cup\overline{(G,A)}\right)^{\mathbb{C}}=\overline{(F,A)}^{\mathbb{C}}\Cap\overline{(G,A)}^{\mathbb{C}}$
{[}From $Proposition\:2.23${]} and this implies that $\overline{(F,A)}\Cup\overline{(G,A)}$
is soft closed set. Therefore, $\overline{(F,A)\Cup(G,A)}\tilde{\subseteq}\overline{(F,A)}\Cup\overline{(G,A)}.$\\
Thus $\overline{(F,A)}\Cup\overline{(G,A)}=\overline{(F,A)\Cup(G,A)}.$
$\newline$\\
$(vii)$ Since $(F,A)\Cap(G,A)\tilde{\subseteq}(F,A)$ and $(F,A)\Cap(G,A)\tilde{\subseteq}(G,A)$,
so by $(v)$, \\
$\overline{(F,A)\Cap(G,A)}\tilde{\subseteq}\overline{(F,A)}$ and
$\overline{(F,A)\Cap(G,A)}\tilde{\subseteq}\overline{(G,A)}$. \\
Thus $\overline{(F,A)\Cap(G,A)}\tilde{\subseteq}\overline{(F,A)}\Cap\overline{(G,A)}$.
$\newline$\\
$\mathbf{Example\:3.15.}$ Here we give an example where $\overline{(F,A)}\Cup\overline{(G,A)}\neq\overline{(F,A)\Cup(G,A)}.$
Consider the soft topological space $(\tilde{X},\tau_{1},A)$ of $Example\:3.12.$
Where the non-null soft closed sets are $(P,A)$ and $(Q,A)$ such
that $P(\alpha)=\{z\},P(\beta)=\{y\}$ and $Q(\alpha)=\{x,y\},Q(\alpha)=\{x\}.$
Then $(P,A)=\overline{(P,A)}$ and $(Q,A)=\overline{(Q,A)}.$ Let
$(V,A)=\overline{(P,A)}\Cup\overline{(Q,A)}$. Then $V(\alpha)=\{x,y,z\},V(\beta)=\{x,y\}.$
But $\overline{(P,A)\Cup(Q,A)}=(\tilde{X},A).$Therefore, $\overline{(P,A)}\Cup\overline{(Q,A)}\neq\overline{(P,A)\Cup(Q,A)}.$$\newline$
\\
$\mathbf{Definition\:3.16.}$ Let $(\tilde{X},\tau,A)$ be a soft
topological space. A soft element $\tilde{x}\tilde{\in}(\tilde{X},A)$
is said to be a limiting soft element of a soft set $(F,A)\in S(\tilde{X})$
if $\forall(G,A)\in\tau$ and for any $\alpha\in A,$ $\tilde{x}(\alpha)\in G(\alpha)$
implies $F(\alpha)\cap[G(\alpha)\setminus\tilde{x}(\alpha)]\neq\phi.$
\\
The soft set formed out of all limiting soft elements of $(F,A)\in S(\tilde{X})$
is called the derived soft set of $(F,A)$ and is denoted by $(F,A)^{\prime}.$
\\
The weak soft closure of a soft set $(F,A)\in S(\tilde{X}),$ denoted
by $\overline{(F,A)}^{w},$ is defined by $\overline{(F,A)}^{w}=(F,A)\Cup(F,A)^{\prime}.$
$\newline$ \\
$\mathbf{Proposition\:3.17}$ If a soft set $(F,A)\in S(\tilde{X})$
is soft closed in a soft topological space $(\tilde{X},\tau,A),$
then $(F,A)$ contains all its limiting soft elements as in $Definition\:3.16.$$\newline$\\
$Proof.$ Let $(F,A)\in S(\tilde{X})$ be soft closed and $\tilde{x}\tilde{\notin}(F,A).$
Since $(F,A)$ is soft closed $(F,A)^{C}\in S(\tilde{X})$ and $(F,A)^{\mathbb{C}}$
is soft open.\\
Now,  $\tilde{x}$ either belongs to $(F,A)^{\mathbb{C}}$ or does
not belongs to $(F,A)^{\mathbb{C}}.$ If $\tilde{x}\tilde{\in}(F,A)^{\mathbb{C}},$
then $\tilde{x}(\lambda)\in[(F,A)^{\mathbb{C}}](\lambda)=X\setminus F(\lambda),\forall\lambda\in A$.
Since $F(\lambda)\cap X\setminus F(\lambda)=\phi,\forall\lambda\in A,$
$\tilde{x}$ can not be limiting soft element of $(F,A).$ If $\tilde{x}\tilde{\notin}(F,A)^{\mathbb{C}},$
then $\tilde{x}(\lambda)\in F(\lambda)$ for some $\lambda\in A$
and $\tilde{x}(\mu)\notin F(\mu)$ for some $\mu(\neq\lambda)\in A.$
Then $\tilde{x}(\mu)\in[(F,A)^{\mathbb{C}}](\mu)$ but $F(\mu)\cap[(F,A)^{\mathbb{C}}](\mu)=\phi$.
So, $\tilde{x}$ can not be limiting soft element of $(F,A)$ as per
$Definition\:3.16.$\\
Therefore $(F,A)$ contains all its limiting soft element.$\newline$\\
$\mathbf{Remark\:3.18.}$ Converse of $Proposition\;3.17$ is not
true.\\
Let $X=\{x,y,z\},$ $A=\{\alpha,\beta\}$ and $\tau=\{(\tilde{\Phi},A),(\tilde{X},A),(F,A),(G,A),(H,A)\}$
where $F(\alpha)=\{y,z\},$ $F(\beta)=\{x\}$; $G(\alpha)=\{y\},$
$G(\beta)=\{y,z\}$; $H(\alpha)=\{y,z\},$ $H(\beta)=\{x,y,z\}.$
Then $\tau$ is a soft topology as per $Definition\:3.1.$ Now, let
$(C,A)$ be any soft set where $C(\alpha)=\{x,z\},$ $C(\beta)=\{y,z\}.$
Let $\tilde{\xi}_{1}(\alpha)=x,$ $\tilde{\xi}_{1}(\beta)=y$; $\tilde{\xi}_{2}(\alpha)=x,$
$\tilde{\xi}_{2}(\beta)=z$. Then $\tilde{\xi}_{1},\tilde{\xi}_{2}$
are the only limiting soft elements of $(C,A)$ as per $Definition\;3.16$
and $\tilde{\xi}_{i}\tilde{\in}(C,A),$ $i=1,2$. But $(C,A)$ is
not closed in $(\tilde{X},\tau,A),$ as $(C,A)^{\mathbb{C}}$ is not
soft open. $\newline$\\
$\mathbf{Definition\:3.19.}$ Let $(\tilde{X},\tau,A)$ be a soft
topological space and $(F,A)\in S(\tilde{X})$. A soft element $\tilde{x}\tilde{\in}(F,A)$
is said to be an interior soft element of $(F,A)$ if $\exists(G,A)\in\tau$
such that $\tilde{x}\tilde{\in}(G,A)\tilde{\subseteq}(F,A).$ $\newline$\\
$\mathbf{Definition\:3.20.}$ The interior of a soft set $(F,A)$
is defined to be the set consisting of all interior soft elements
of $(F,A).$ The interior of the soft set $(F,A)$ denoted by $Int(F,A).$
Thus,\\
 $Int(F,A)=\{\tilde{x}\tilde{\in}(F,A):\tilde{x}\tilde{\in}(G,A)\tilde{\subseteq}(F,A)\: for\: some\:(G,A)\in\tau\}.$\\
 $SS[Int(F,A)]$ is said to be soft interior of $(F,A)$ and denoted
by $(F,A)^{\textdegree}.$ $\newline$\\
$\mathbf{Proposition\:3.21.}$ Let $(\tilde{X},\tau,A)$ be a soft
topological space and $(F,A),(G,A)\in S(\tilde{X})$. Then\\
$(i)$ $(F,A)^{\textdegree}\tilde{\subseteq}(F,A).$\\
$(ii)$ $(F,A)\tilde{\subseteq}(G,A)\Rightarrow(F,A)^{\textdegree}\tilde{\subseteq}(G,A)^{\textdegree}.$\\
$(iii)$ $Int[(F,A)\Cup(G,A)]\supseteq Int(F,A)\cup Int(G,A)$. Also,
$[(F,A)\Cup(G,A)]^{\textdegree}\tilde{\supseteq}(F,A)^{\textdegree}\Cup(G,A)^{\textdegree}.$\\
$(iv)$ $Int(F,A)\cap Int(G,A)\subseteq Int[(F,A)\Cap(G,A)].$$\newline$\\
$\mathbf{Proposition\:3.22.}$ Let $(\tilde{X},\tau,A)$ be a soft
topological space. Then $(F,A)(\neq(\tilde{\Phi},A))\in S(\tilde{X})$
is soft open if and only if every $\tilde{x}\tilde{\in}(F,A)$ is
an interior soft element of $(F,A).$ $\newline$\\
$\mathbf{Proposition\:3.23.}$ Let $(\tilde{X},\tau,A)$ be a soft
topological space and $(F,A)\in S(\tilde{X})$. Then $(F,A)^{\textdegree}$
is the elementary union of all soft open sets contained in $(F,A).$
It is the largest soft open set in $(\tilde{X},\tau,A)$ contained
in $(F,A).$ $\newline$\\
$\mathbf{Definition\:3.24.}$ Let $(\tilde{X},\tau,A)$ be a soft
topological space. Then $(F,A)(\neq(\tilde{\Phi},A))\in S(\tilde{X})$
is a soft neighbourhood (soft nbd) of the soft element $\tilde{x}$
if there exists a soft set $(G,A)\in\tau$ such that $\tilde{x}\tilde{\in}(G,A)\tilde{\subseteq}(F,A).$
The soft nbd system at a soft element $\tilde{x}$, denoted by $\aleph_{\tau}(\tilde{x})$,
is the family of all its soft nbds. $\newline$\\
$\mathbf{Proposition\:3.25.}$ Let $(\tilde{X},\tau,A)$ be a soft
topological space. Then $(F,A)(\neq(\tilde{\Phi},A))\in S(\tilde{X})$
is soft open if and only if $(F,A)$ is soft nbd of all of its soft
elements.$\newline$\\
$\mathbf{Proposition\:3.26.}$ Let $(\tilde{X},\tau,A)$ be a soft
topological space and for $\tilde{x}\in SE(\tilde{X})$, let $\aleph_{\tau}(\tilde{x})$
be the soft nbd system at the soft element $\tilde{x}$. Then\\
 $(i)$ $\aleph_{\tau}(\tilde{x})\neq\phi,$ $\forall\tilde{x}\in SE(\tilde{X})$

\noindent $(ii)$ $\tilde{x}\tilde{\in}(F,A),\forall(F,A)\in\aleph_{\tau}(\tilde{x})$ 

\noindent $(iii)$ $(F,A)\in\aleph_{\tau}(\tilde{x}),(F,A)\tilde{\subseteq}(G,A)\Rightarrow(G,A)\in\aleph_{\tau}(\tilde{x})$ 

\noindent $(iv)$ $(F,A),(G,A)\in\aleph_{\tau}(\tilde{x})\Rightarrow(F,A)\Cap(G,A)\in\aleph_{\tau}(\tilde{x})$

\noindent $(v)$ $(F,A)\in\aleph_{\tau}(\tilde{x})\Rightarrow\exists(G,A)\in\aleph_{\tau}(\tilde{x})$
such that $(G,A)\tilde{\subseteq}(F,A)$ and $(G,A)\in\aleph_{\tau}(\tilde{y}),\forall\tilde{y}\tilde{\in}(G,A)$.$\newline$\\
$\mathbf{Definition\:3.27.}$ A mapping $\nu:SE(\tilde{X})\rightarrow P(S(\tilde{X}))$
is said to be a soft nbd operator on $SE(\tilde{X})$ if the following
conditions hold:\\
$(N1)$ $\nu(\tilde{x})\neq\phi,$ $\forall\tilde{x}\in SE(\tilde{X})$

\noindent $(N2)$ $\tilde{x}\tilde{\in}(F,A),\forall(F,A)\in\nu(\tilde{x})$ 

\noindent $(N3)$ $(F,A)\in\nu(\tilde{x}),(F,A)\tilde{\subseteq}(G,A)\Rightarrow(G,A)\in\nu(\tilde{x})$ 

\noindent $(N4)$ $(F,A),(G,A)\in\nu(\tilde{x})\Rightarrow(F,A)\Cap(G,A)\in\nu(\tilde{x})$

\noindent $(N5)$ $(F,A)\in\nu(\tilde{x})\Rightarrow\exists(G,A)\in\nu(\tilde{x})$
such that $(G,A)\tilde{\subseteq}(F,A)$ and $(G,A)\in\nu(\tilde{y}),\forall\tilde{y}\tilde{\in}(G,A)$.$\newline$\\
$\mathbf{Example\:3.28.}$ If $(\tilde{X},\tau,A)$ is a soft topological
space, then the mapping $\nu:SE(\tilde{X})\rightarrow P(S(\tilde{X}))$
defined by $\nu(\tilde{x})=\aleph_{\tau}(\tilde{x})$, where $\aleph_{\tau}(\tilde{x})$
is the soft nbd system at the soft element $\tilde{x}$, is a soft
nbd operator on $SE(\tilde{X})$. $\newline$\\
\section{Soft base of a soft topology }
$\mathbf{Definition\:4.1.}$ Let $(\tilde{X},\tau,A)$ be a soft topological
space. Then a subcollection $\text{\ensuremath{\mathscr{B}}}$ of
$\tau,$ containing $(\tilde{\Phi},A)$, is said to be an open base
$\tau$ iff $\forall\tilde{x}\tilde{\in}(\tilde{X},A)$ and for any
soft open set $(F,A)$ containing the soft element $\tilde{x}$, there
exists $(G,A)\in\text{\ensuremath{\mathscr{B}}}$ such that $\tilde{x}\tilde{\in}(G,A)\tilde{\subseteq}(F,A)$.$\newline$\\
$\mathbf{Examples\:4.2.}$ Let $X=\{x,y,z,t\},$ $A=\{\alpha,\beta\}$
and $\tau=\{(\tilde{\Phi},A),(\tilde{X},A),$ $(F_{1},A),(F_{2},A),(F_{3},A),(F_{4},A),(F_{5},A),(F_{6},A),(F_{7},A)\}$
where\\
$F_{1}(\alpha)=\{x\},$ $F_{1}(\beta)=\{t\}$; $F_{2}(\alpha)=\{y\},$
$F_{2}(\beta)=\{z\}$; \\
$F_{3}(\alpha)=\{t\},$ $F_{3}(\beta)=\{x\}$; $F_{4}(\alpha)=\{x,y\},$
$F_{4}(\beta)=\{z,t\}$; \\
$F_{5}(\alpha)=\{y,t\},$ $F_{5}(\beta)=\{z,x\};$ $F_{6}(\alpha)=\{x,t\},$
$F_{6}(\beta)=\{x,t\};$\\
$F_{7}(\alpha)=\{x,y,t\},$ $F_{7}(\beta)=\{x,z,t\}$. \\
Then $\tau$ is a soft topology on $(\tilde{X},A)$. \\
Let $\text{\ensuremath{\mathscr{B}}}=\{(\tilde{\Phi},A),(\tilde{X},A),(F_{1},A),(F_{2},A),(F_{3},A),(F_{4},A),(F_{5},A),(F_{6},A)\}.$
\\
Then $\text{\ensuremath{\mathscr{B}}}$ forms an open base for $\tau.$
$\newline$\\
$\mathbf{Proposition\:4.3.}$ Let $(\tilde{X},\tau,A)$ be a soft
topological space and $\text{\ensuremath{\mathscr{B}}}$ is an open
base for $\tau$. Then every member of $\tau$ can be expressed as
the elementary union of some members of $\text{\ensuremath{\mathscr{B}}}.$
$\newline$\\
$\mathbf{Remark\:4.4.}$ Converse of $Proposition\:4.3$ is not true.
Consider the soft topological space $(\tilde{X},\tau,A)$ of $Example\;4.2$
and let $\text{\ensuremath{\mathscr{B}}}=\{(\tilde{\Phi},A),(\tilde{X},A),(F_{1},A),$
$(F_{2},A),(F_{3},A)\}.$ \\
Then $\text{\ensuremath{\mathscr{B}}}$ satisfies the condition of
$Proposition\:4.3$ but $\text{\ensuremath{\mathscr{B}}}$ is not
an open base for $\tau$ as the soft element $\bar{x}\tilde{\in}(F_{6},A)$
but there is no soft set in $\text{\ensuremath{\mathscr{B}}}$ containing
$\bar{x}$ and contained in $(F_{6},A).$ $\newline$\\
$\mathbf{Proposition\:4.5.}$ If a collection $\text{\ensuremath{\mathscr{B}}}$
of soft sets of $S(\tilde{X})$ forms an open base of a soft topological
space $(\tilde{X},\tau,A)$, then the following conditions are satisfied:\\
$(i)$ $(\tilde{\Phi},A)\in\text{\ensuremath{\mathscr{B}}}.$\\
$(ii)$ $(\tilde{X},A)$ is elementary union of some members of $\text{\ensuremath{\mathscr{B}}}.$\\
$(iii)$ If $(F_{1},A),(F_{2},A)\in\text{\ensuremath{\mathscr{B}}}$
and $\tilde{x}\tilde{\in}(F_{1},A)\Cap(F_{2},A)$, then there exists
$(F_{3},A)\in\text{\ensuremath{\mathscr{B}}}$ such that $\tilde{x}\tilde{\in}(F_{3},A)\tilde{\subseteq}(F_{1},A)\Cap(F_{2},A)$.$\newline$\\
$\mathbf{Remark\:4.6.}$ Converse of $Proposition\:4.5$ is not true.
Consider the soft topological space $(\tilde{X},\tau,A)$ of $Example\:4.2$
and let $\text{\ensuremath{\mathscr{B}}}$ as in $Remark\:4.4.$ Then
$\text{\ensuremath{\mathscr{B}}}$ satisfies all the condition of
$Proposition\:4.5$ but $\text{\ensuremath{\mathscr{B}}}$ is not
base for the soft topology $\tau$.$\newline$\\
$\mathbf{Definition\:4.7.}$ A family $\mathcal{S}$ of subsets of
$X$ is said to be a sub-base for a soft topological space $(\tilde{X},\tau,A)$,
if the family of all finite elementary intersection of members of
$\mathcal{S}$ is a soft base for $\tau.$ $\newline$
\section{Soft function and soft continuous function}
Proceeding as in \cite{key-23}, where definition of soft mapping
has been given using 'soft point' concept, we introduce here a definition
of soft function using the concept of 'soft element'.$\newline$

\noindent $\mathbf{Definition\:5.1.}$ Let $X$ and $Y$ be two non-empty
sets and $\{f_{\lambda}:X\rightarrow Y,\lambda\in A\}$ be a collection
of functions. Then a function $f:SE(\tilde{X})\rightarrow SE(\tilde{Y})$
defined by $[f(\tilde{x})](\lambda)=f_{\lambda}(\tilde{x}(\lambda)),\forall\lambda\in A$
is called a soft function. $\newline$\\
$\mathbf{Definition\:5.2.}$ Let $f:SE(\tilde{X})\rightarrow SE(\tilde{Y})$
be a soft function. Then \\
$(i)$ image of a soft set $(F,A)$ over $X$ under the soft function
$f$, denoted by $f[(F,A)]$, is defined by $f[(F,A)]=SS\{f(SE(F,A))\}$
i.e. $f[(F,A)](\lambda)$ $=f_{\lambda}(F(\lambda)),$ $\forall\lambda\in A.$\\
$(ii)$ inverse image of a soft set $(G,A)$ over $Y$ under the soft
function $f$, denoted by $f^{-1}[(G,A)]$, is defined by $f^{-1}[(G,A)]=SS\{f^{-1}(SE(G,A))\}$
i.e. $f^{-1}[(G,A)](\lambda)=f_{\lambda}^{-1}(G(\lambda)),\forall\lambda\in A.$
$\newline$\\
$\mathbf{Definition\:5.3.}$ Let $f:SE(\tilde{X})\rightarrow SE(\tilde{Y})$
be a soft function associated with the family of functions $\{f_{\lambda}:X\rightarrow Y,\lambda\in A\}$.
Then $f$ is said to be
\noindent $(i)$ injective if $\tilde{x}\neq\tilde{y}$ implies $f(\tilde{x})\neq f(\tilde{y}).$
\\
$(ii)$ surjective if $f(\tilde{X},A)=(\tilde{Y},A).$\\
$(iii)$ bijective if both injective and surjective.$\newline$\\
$\mathbf{Proposition\:5.4.}$ Let $X$ and $Y$ be two non-empty sets
and $A$ be the parameter set. Also let $f:SE(\tilde{X})\rightarrow SE(\tilde{Y})$
be a soft function associated with the family of functions $\{f_{\lambda}:X\rightarrow Y,\lambda\in A\}.$
If $(F,A)\in S(\tilde{X})$ then \\
$(i)$ $ff^{-1}(F,A)\tilde{\subseteq}(F,A)$.\\
$(ii)$ $(F,A)\tilde{\subseteq}f^{-1}f(F,A)$.$\newline$\\
$\mathbf{Proposition\:5.5.}$ Let $X$ and $Y$ be two non-empty sets
and $A$ be the parameter set. Also let $f:SE(\tilde{X})\rightarrow SE(\tilde{Y})$
be a soft function associated with the family of functions $\{f_{\lambda}:X\rightarrow Y,\lambda\in A\}.$
If $(F_{1},A),(F_{2},A)\in S(\tilde{X})$ then \\
$(i)$ $(F_{1},A)\tilde{\subseteq}(F_{2},A)\Rightarrow f[(F_{1},A)]\tilde{\subseteq}f[(F_{2},A)].$\\
$(ii)$ $f[(F_{1},A)\Cup(F_{2},A)]=f[(F_{1},A)]\Cup f[(F_{2},A)]$.\\
$(iii)$ $f[(F_{1},A)\Cap(F_{2},A)]\tilde{\subseteq}f[(F_{1},A)]\Cap f[(F_{2},A)]$.\\
$(iv)$ $f[(F_{1},A)\Cap(F_{2},A)]=f[(F_{1},A)]\Cap f[(F_{2},A)]$,
if $f$ is one-one.$\newline$\\
$\mathbf{Proposition\:5.6.}$ Let $X$ and $Y$ be two non-empty sets
and $A$ be the parameter set. Also let $f:SE(\tilde{X})\rightarrow SE(\tilde{Y})$
be a soft function associated with the family of functions $\{f_{\lambda}:X\rightarrow Y,\lambda\in A\}.$
If $(F_{1},A),(F_{2},A)\in S(\tilde{X})$ then \\
$(i)$ $(F_{1},A)\tilde{\subseteq}(F_{2},A)\Rightarrow f^{-1}[(F_{1},A)]\tilde{\subseteq}f^{-1}[(F_{2},A)].$\\
$(ii)$ $f^{-1}[(F_{1},A)\Cup(F_{2},A)]=f^{-1}[(F_{1},A)]\Cup f^{-1}[(F_{2},A)]$.\\
$(iii)$ $f^{-1}[(F_{1},A)\Cap(F_{2},A)]=f^{-1}[(F_{1},A)]\Cap f^{-1}[(F_{2},A)]$.
$\newline$\\
$\mathbf{Definition\:5.7.}$ Let $(\tilde{X},\tau,A)$ and $(\tilde{Y},\nu,A)$
be two soft topological spaces and $f:SE(\tilde{X})\rightarrow SE(\tilde{Y})$
be a soft function associated with the family of functions $\{f_{\lambda}:X\rightarrow Y,\lambda\in A\}.$
Then we denote this soft function as $f:(\tilde{X},\tau,A)\rightarrow(\tilde{Y},\nu,A)$.
\\
Now $f:(\tilde{X},\tau,A)\rightarrow(\tilde{Y},\nu,A)$ is said to
be soft continuous at $\tilde{x}_{0}\tilde{\in}(\tilde{X,}A),$ if
for every $(V,A)\in\nu$ such that $f(\tilde{x}_{0})\tilde{\in}(V,A),$
there exists $(U,A)\in\tau$ such that $\tilde{x}_{0}\tilde{\in}(U,A)$
and $f(U,A)\tilde{\subseteq}(V,A).$\\
$f$ is said to be soft continuous on $(\tilde{X},\tau,A)$ if it
is soft continuous at each soft element $\tilde{x}_{0}\tilde{\in}(\tilde{X,}A).$
$\newline$\\
$\mathbf{Proposition\:5.8.}$ Let $(\tilde{X},\tau,A)$ and $(\tilde{Y},\nu,A)$
be two soft topological spaces and $f:(\tilde{X},\tau,A)\rightarrow(\tilde{Y},\nu,A)$
be a soft function. Then the followings are related as follows:\\
$(i)\Leftrightarrow(ii),$ $(ii)\Leftrightarrow(iii)$ and $(ii)\Rightarrow(iv).$
\\
$(i)$ $f$ is soft continuous.\\
$(ii)$ For all $(V,A)\in\nu,$ $f^{-1}(V,A)\in\tau.$\\
$(iii)$ There exists a subbase $\wp$ for $\nu$ such that $f^{-1}(V,A)\in\tau$
for all $(V,A)\in\wp.$\\
$(iv)$ for any closed soft set $(F,A)\in S(\tilde{Y})$ in $(\tilde{Y},\nu,A)$,
$f^{-1}(F,A)$ is soft closed in $(\tilde{X},\tau,A)$ . $\newline$\\
$\mathbf{Remark\:5.9.}$ In $Proposition\:5.8$ $(iv)\Rightarrow(i)$
does not hold.\\
Let $X=\{x,y,z\}$ and $A=\{\alpha,\beta\}.$ Let $\tau_{1}=\{(\tilde{\Phi},A),(\tilde{X},A),(F,A)\}$
and $\tau_{2}=\{(\tilde{\Phi},A),(\tilde{X},A)\},$ where $F(\alpha)=\{x,y,z\},F(\beta)=\{x,z\}.$
Then $\tau_{1}$ and $\tau_{2}$ are soft topologies on $(\tilde{X},A).$
Consider the soft function $i:(\tilde{X},\tau_{2},A)\rightarrow(\tilde{X},\tau_{1},A)$
corresponding to the identity function $i:X\rightarrow X$. Then $i^{-1}:(\tilde{X},\tau_{1},A)\rightarrow(\tilde{X},\tau_{2},A)$
maps all soft closed sets of $\tau_{1}$ to soft closed sets of $\tau_{2}.$
But $i^{-1}(F,A)\notin\tau_{2}.$ Therefore, $i:(\tilde{X},\tau_{2},A)\rightarrow(\tilde{X},\tau_{1},A)$
is not soft continuous function.$\newline$\\
$\mathbf{Definition\:5.10.}$ A soft function $f:(\tilde{X},\tau,A)\rightarrow(\tilde{Y},\nu,A)$
is said to be \\
$(i)$ soft open if $f$ maps soft open sets of $\tau$ to soft open
sets of $\nu.$\\
$(ii)$ soft closed if $f$ maps soft closed sets of $\tau$ to soft
closed sets of $\nu.$ $\newline$ \\
$\mathbf{Definition\:5.11.}$ A soft function $f:(\tilde{X},\tau,A)\rightarrow(\tilde{Y},\nu,A)$
is said to be soft homeomorphism if\\
$(i)$ $f$ is bijective\\
$(ii)$ $f$ is soft continuous\\
$(iii)$ $f^{-1}$ is soft continuous. $\newline$\\
$\mathbf{Proposition\:5.12.}$ Let $f:(\tilde{X},\tau,A)\rightarrow(\tilde{Y},\nu,A)$
be a soft function. Then the followings are equivalent:\\
$(i)$ $f$ is a soft homeomorphism.\\
$(ii)$ $f$ is bijective, $f$ , $f^{-1}$ are soft continuous. \\
$(iii)$ $f$ is bijective, soft open and soft continuous.\\
$(iv)$ $f^{-1}$ is a soft homeomorphism.$\newline$\\
\section{Soft separation axioms}
$\mathbf{Definition\:6.1.}$ Let $(\tilde{X},\tau,A)$ be a soft topological
space. If for $\tilde{x},\tilde{y}\in SE(\tilde{X})$ with $\tilde{x}(\lambda)\neq\tilde{y}(\lambda),\forall\lambda\in A$,\\
 $\text{(\mbox{i})}$ there exists $(F,A)\in\tau$ such that {[}$\tilde{x}(\lambda)\in(F,A)(\lambda)$
and $\tilde{y}(\lambda)\notin(F,A)(\lambda)${]} or {[}$\tilde{y}(\lambda)\in(F,A)(\lambda)$
and $\tilde{x}(\lambda)\notin(F,A)(\lambda)${]}, $\forall\lambda\in A$,
then $(\tilde{X},\tau,A)$ is called a soft $T_{0}$ space.$\newline$\\
$\text{(\mbox{ii})}$ there exist $(F,A),(G,A)\in\tau$ such that
{[}$\tilde{x}(\lambda)\in(F,A)(\lambda)$, $\tilde{y}(\lambda)\notin(F,A)(\lambda)${]}
and {[}$\tilde{y}(\lambda)\in(G,A)(\lambda)$, $\tilde{x}\notin(G,A)(\lambda)${]},
$\forall\lambda\in A$ then $(\tilde{X},\tau,A)$ is called a soft
$T_{1}$ space.$\newline$\\
$\text{(\mbox{iii})}$ there exist $(F,A),(G,A)\in\tau$ such that
$\tilde{x}\tilde{\in}(F,A)$, $\tilde{y}\tilde{\in}(G,A)$ and $(F,A)\tilde{\cap}(G,A)$
~$=(\tilde{\Phi},A),$ then $(\tilde{X},\tau,A)$ is called a soft
$T_{2}$ space.$\newline$\\
$\mathbf{Remark\:6.2.}$ $(\text{\mbox{i}})$ Every soft $T_{1}$
space is a soft $T_{0}$ space.\\
$\text{(\mbox{ii})}$ Every soft $T_{2}$ space is a soft $T_{1}$
space.$\newline$\\
$\mathbf{Example\:6.3.}$ Let $\mathbb{R}$ be the real number space
and $A$ be a non-empty parameter set. Let for each $\alpha\in A,$
$\tau(\alpha)$ be the usual topology on $\mathbb{R}.$ Then $(\tilde{\mathbb{R}},\tau,A)$
where $\tau$ be the soft topology generated by $\tau(\alpha)$ as
in $Proposition\;3.9,$ is a soft $T_{2}$ space. $\newline$ \\
$\mathbf{Definition\:6.4.}$ If a soft set of $S(\tilde{X})$ contains
exactly one soft element $\tilde{x,}$ then we denote this soft set
by $(\tilde{x},A)$ i.e. $(\tilde{x},A)(\lambda)=\tilde{x}(\lambda),\forall\lambda\in A.$
$\newline$ \\
$\mathbf{Proposition\:6.5.}$ Let $(\tilde{X},\tau,A)$ be a soft
$T_{1}$ topological space. Then for any $\tilde{x}\in SE(\tilde{X})$,
$(\tilde{x},A)$ is soft closed. $\newline$ \\
$Proof.$ Let $\tilde{x}\in SE(\tilde{X})$. We claim that $(\tilde{x},A)^{C}=(\tilde{x},A)^{\mathbb{C}}$
is a soft nbd of each of its soft element. For, let $\tilde{y}\tilde{\in}(\tilde{x},A)^{C}.$
Then $\tilde{x}(\lambda)\neq\tilde{y}(\lambda),\forall\lambda\in A$
and since $(\tilde{X},\tau,A)$ is soft $T_{1}$, there exists $(G,A)\in\tau$
such that {[}$\tilde{y}(\lambda)\in(G,A)(\lambda)$, $\tilde{x}\notin(G,A)(\lambda)${]},
$\forall\lambda\in A$. But this means $(G,A)\tilde{\subseteq}(\tilde{x},A)^{C}$
and hence $(\tilde{x},A)^{C}$ is a soft nbd of $\tilde{y}.$ So,
$(\tilde{x},A)^{C}=(\tilde{x},A)^{\mathbb{C}}$ is open. Hence $(\tilde{x},A)$
is soft closed. $\newline$ \\
$\mathbf{Definition\:6.6.}$ A soft topological space $(\tilde{X},\tau,A)$
is said to be a soft regular space if for any soft closed set $(F,A)$
and any soft element $\tilde{x}$ such that $\tilde{x}(\lambda)\notin(F,A)(\lambda)$,
$\forall\lambda\in A,$ $\exists$ $(G,A),(H,A)\in\tau$ such that
$(F,A)\tilde{\subseteq}(G,A),\tilde{x}\tilde{\in}(H,A)$ and $(F,A)\Cap(G,A)=(\tilde{\Phi},A).$
\\
If in addition,$(\tilde{X},\tau,A)$ is soft $T_{1},$ then $(\tilde{X},\tau,A)$
is called soft $T_{3}$ space. $\newline$\\
$\mathbf{Example\:6.7.}$ Let $X=\{x,y,z\},A=\{\alpha,\beta\}$ and
$\tau=\{(\tilde{\Phi},A),(\tilde{X},A),(F,A),$ $(G,A)\},$ where
$F(\alpha)=\{x,z\},F(\beta)=\{y\};$ $G(\alpha)=\{y\},G(\beta)=\{x,z\}.$
Then $(\tilde{X},\tau,A)$ is a soft regular space. $\newline$\\
$\mathbf{Proposition\:6.8.}$ Let $(\tilde{X},\tau,A)$ be a soft
topological space. If for all $\tilde{x}\in SE(\tilde{X})$ and for
all soft open sets $(U,A)$ such that $\tilde{x}\tilde{\in}(U,A),$
$\exists(V,A)\in\tau$ such that $\tilde{x}\tilde{\in}(V,A)\tilde{\subseteq}\overline{(V,A)}\tilde{\subseteq}(U,A)$,
then $(\tilde{X},\tau,A)$ is soft regular, where closure is taken
as per $Definition\:3.13$. $\newline$\\
$Proof.$ Let $\tilde{x}\in SE(\tilde{X})$ and $(F,A)$ be any soft
closed set such that $\tilde{x}(\lambda)\notin(F,A)(\lambda)$, $\forall\lambda\in A.$
Consider $(U,A)=(F,A)^{C}.$ Then $(U,A)\in\tau,$ as $(U,A)\in S(\tilde{X})$
and $(U,A)=(F,A)^{C}=(F,A)^{\mathbb{C}}$ and $\tilde{x}\tilde{\in}(U,A).$
So by given condition $\exists(V,A)\in\tau$ such that $\tilde{x}\tilde{\in}(V,A)\tilde{\subseteq}\overline{(V,A)}\tilde{\subseteq}(U,A)$.
Let $(W,A)=\overline{(V,A)}^{C}.$ Then $(W,A)\in\tau,$ as $(W,A)\in S(\tilde{X})$
and $(W,A)=\overline{(V,A)}^{C}=\overline{(V,A)}^{\mathbb{C}}.$ Thus
$(F,A)=(U,A)^{C}$ $\tilde{\subseteq}(\tilde{X},A)\setminus\overline{(V,A)}=(W,A).$
Also, $(V,A)\tilde{\cap}(W,A)=(\tilde{\Phi},A)$ i.e. $(V,A)\Cap(W,A)=(\tilde{\Phi},A).$
Therefore, $(\tilde{X},\tau,A)$ is a soft regular space. $\newline$
\\
 The following example shows that the converse of $Proposition\:6.8$
is not true.$\newline$ \\
$\mathbf{Example\:6.9.}$ Let $X=\{x,y\},A=\{\alpha,\beta\}$ and
$\tau=\{(\tilde{\Phi},A),(\tilde{X},A),(F_{1},A),$ $(F_{2},A),(F_{3},A),(F_{4},A)\},$
where $F_{1}(\alpha)=\{x,y\},F_{1}(\beta)=\{x\};$ $F_{2}(\alpha)=\{y\},$
$F_{2}(\beta)=\{x,y\};$ $F_{3}(\alpha)=\{y\},F_{3}(\beta)=\{x\};$
$F_{4}(\alpha)=\{x\},F_{4}(\beta)=\{y\}.$ Then $(\tilde{X},\tau,A)$
is a soft regular space. But $(F_{1},A)$ is a soft open set containing
the soft element $\overline{x}$ for which the given condition of
$Proposition\:6.8$ is not satisfied. $\newline$\\
$\mathbf{Definition\:6.10.}$ A soft topological space $(\tilde{X},\tau,A)$
is said to be a soft normal space if for any two soft closed sets
$(F,A)$ and $(G,A)$ such that $(F,A)\tilde{\cap}(G,A)=(\tilde{\Phi},A)$,
$\exists$ $(U,A),(V,A)\in\tau$ such that $(F,A)\tilde{\subseteq}(U,A),(G,A)\tilde{\subseteq}(V,A)$
and $(U,A)\Cap(V,A)=(\tilde{\Phi},A).$ \\
If in addition,$(\tilde{X},\tau,A)$ is soft $T_{1},$ then $(\tilde{X},\tau,A)$
is called soft $T_{4}$ space.$\newline$\\
$\mathbf{Proposition\:6.11.}$ Let $(\tilde{X},\tau,A)$ be a soft
topological space. If for all soft closed sets $(F,A)$ and for all
soft open sets $(U,A)$ such that $(F,A)\tilde{\subseteq}(U,A),$
$\exists(V,A)\in\tau$ such that $(F,A)\tilde{\subseteq}(V,A)\tilde{\subseteq}\overline{(V,A)}\tilde{\subseteq}(U,A)$,
then $(\tilde{X},\tau,A)$ is soft normal, where closure is taken
as per $Definition\:3.13$. $\newline$\\
$Proof.$ Let the given condition be satisfied.\\
Let $(F_{1},A),(F_{2},A)$ be two soft closed sets such that $(F_{1},A)\tilde{\cap}(F_{2},A)=(\tilde{\Phi},A)$.
Consider $(U,A)=(F_{2},A)^{C}$. Then $(U,A)\in\tau,$ as $(U,A)\in S(\tilde{X})$
and $(U,A)=(F_{2},A)^{C}=(F_{2},A)^{\mathbb{C}}.$ Also, $(F_{1},A)\tilde{\subseteq}(U,A).$
So, by given condition $\exists(V,A)\in\tau$ such that $(F_{1},A)\tilde{\subseteq}(V,A)\tilde{\subseteq}\overline{(V,A)}\tilde{\subseteq}(U,A)$.
Let $(W,A)=\overline{(V,A)}^{C}.$ Then $(W,A)\in\tau,$ as $(W,A)\in S(\tilde{X})$
and $(W,A)=\overline{(V,A)}^{C}=\overline{(V,A)}^{\mathbb{C}}.$ Thus
$(F_{2},A)=(\tilde{X},A)\setminus(U,A)$ $\tilde{\subseteq}(\tilde{X},A)\setminus\overline{(V,A)}=(W,A).$
Also, $(V,A)\tilde{\cap}(W,A)=(\tilde{\Phi},A)$ i.e. $(V,A)\Cap(W,A)=(\tilde{\Phi},A).$
Therefore, $(\tilde{X},\tau,A)$ is a soft normal space. $\newline$\\
The following example shows that the converse of $Proposition\:6.11$
is not true. $\newline$ \\
$\mathbf{Example\:6.12.}$ Let $X=\{x,y\},A=\{\alpha,\beta\}$ and
$\tau=\{(\tilde{\Phi},A),(\tilde{X},A),(F_{1},A),$ $(F_{2},A),\},$where
$F_{1}(\alpha)=\{x,y\},F_{1}(\beta)=\{x\};$ $F_{2}(\alpha)=\{x\},F_{2}(\beta)=\{y\}.$
Then $(\tilde{X},\tau,A)$ is a soft normal space. But $(F_{1},A)$
is a soft open set containing the soft closed set $(F_{3},A)$, (say),
where $F_{3}(\alpha)=\{y\},F_{3}(\beta)=\{x\}$ for which the given
condition of $Proposition\;6.11$ is not satisfied. $\newline$\\

\end{document}